\documentclass[12pt]{article}
\usepackage[english]{babel}
\usepackage{amsmath,amsfonts,amssymb,amsthm} % math, align* environment

\usepackage{mathdots}

\usepackage{epsf}
\usepackage{graphicx} 

\usepackage[numbers]{natbib}

\usepackage{subfigure}
\usepackage{float}
\usepackage{tikz}
\usetikzlibrary{positioning, calc, fit, decorations.pathreplacing}

\newtheorem{remark}{Remark}

\usepackage{url}
\usepackage{xcolor}
\usepackage{hyperref} % to have hyperlinks in the pdf document
\newcommand\restr[2]{{% we make the whole thing an ordinary symbol
  \left.\kern-\nulldelimiterspace % automatically resize the bar with \right
  #1 % the function
  \vphantom{\big|} % pretend it's a little taller at normal size
  \right|_{#2} % this is the delimiter
  }} % use as "\restr{f}{A}"

\title{On the recent advances of spectral analysis for systems arising from fully-implicit RK methods}
\author{M. Outrata}
\date{\vspace{-5ex}}

\begin{document}
\maketitle

\begin{abstract}
This work deals with two groups of spectral analysis results for matrices arising in fully implicit Runge-Kutta methods used for linear time-dependent partial differential equations. These were applied for different formulations of the same problem and used different tools to arrive at results that do not immediately coincide. We show the equivalence of the results as well as the equivalence of the approaches, unifying the two directions.
\end{abstract}

\section{Introduction}\label{sec_Intro}
In recent years we have seen a renewed interest in using fully implicit Runge-Kutta (IRK) methods as numerical integrators for (linear) time-dependent partial differential equations (PDEs), see~\cite{Howle2021new,farrell2021irksome,southworth2022fastI,southworth2022fastII,
dravins2024spectral,pearson2024pint, outrata2023irksstage},
among others and also~\cite{wanner1987solving,wanner1996solving} for a more complete overview of Runge-Kutta methods in the context of numerical solvers for ordinary differential equations (ODEs). We will consider a time-dependent PDE in the form

\begin{equation}\label{eqn_secIntro_ParabPDE}
\begin{gathered}
\frac{\partial}{\partial t} u = \mathcal{L}u + f \quad \mathrm{in} \; \Omega\times (0,T),\\
\mathcal{B} u(\mathbf{x},t) = g(\mathbf{x},t) \quad \mathrm{on} \;  \partial\Omega \times (0,T),
\qquad
u(\mathbf{x},0) = u^{(0)}(\mathbf{x}) \quad \mathrm{in} \; \Omega,
\end{gathered}
\end{equation}

\noindent on a given connected, bounded domain $\Omega \subset \mathbb{R}^d$ and $(0,T)\subset \mathbb{R}$ for some given data $f, g$ and $u^{(0)}$. Moreover, we will assume that the spatial operator $\mathcal{L}$ is bounded, coercive and self-adjoint. We note that these properties can be relaxed to other relevant cases, see~\cite{outrata2023irksstage,dravins2024spectral}. We first discretize in space using, e.g., a finite elements (FE) scheme, obtaining a system of $n$ ODEs

\begin{equation*}
\frac{\partial}{\partial t} M \mathbf{u}(t) = K \mathbf{u}(t) + \mathbf{b}^{(ST)}(t),
\quad \mathrm{with} \quad
\mathbf{u}(0) = \mathbf{u}^{(0)},
\end{equation*}

\noindent where $M,K \in \mathbb{R}^{n\times n}$ are the mass and stiffness matrices (corresponding to the chosen FE scheme and the operator $\mathcal{L}$), the vector function $\mathbf{b}^{(ST)}(t)$ aggregates the contributions of both $f(t)$ and $g(t)$ and the vector $\mathbf{u}^{(0)}$ corresponds to the initial condition $u^{(0)}(\mathbf{x})$. Any IRK method is then given by the number of stages $s\in \mathbb{N}$ and its Butcher table, succinctly written as

\begin{table}[h]
\centering
 \begin{tabular}{l|l}
 $\mathbf{c}$	& $A$ \\
\hline
 & $\mathbf{b}^T$~~\\[1pt]
\end{tabular}, $\quad A \in \mathbb{R}^{s\times s}, \quad \mathrm{and} \quad \mathbf{b},\mathbf{c} \in \mathbb{R}^s$.
\end{table}

\noindent For a given timestep $\tau$, the IRK method progresses the solution forward in time at timepoints $t_m := \tau m$ using the approximation $\mathbf{u}(t_m) \approx \mathbf{u}^{(m)}$ with

\begin{equation*}
\mathbf{u}^{(m)} := \mathbf{u}^{(m-1)} + \tau \sum\limits_{i=1}^{s} \mathbf{b}_i \mathbf{k}^{(m)}_i,
\end{equation*}

\noindent where the so-called stage-functions $\mathbf{k}^{(m)}_i$ satisfy

\begin{equation}\label{eqn_secIntro_StageEqntsInBasicForm}
M\mathbf{k}^{(m)}_i =\mathbf{b}^{(ST)}(t_{m-1} + c_i\tau) +  K\mathbf{k}^{(m-1)}_i + \tau \sum\limits_{j=1}^{s} a_{ij} K \mathbf{k}^{(m)}_j, \qquad i=1,\dotsc ,s.
\end{equation}

\begin{figure}[t]
\centering
\includegraphics[width=.8\linewidth]{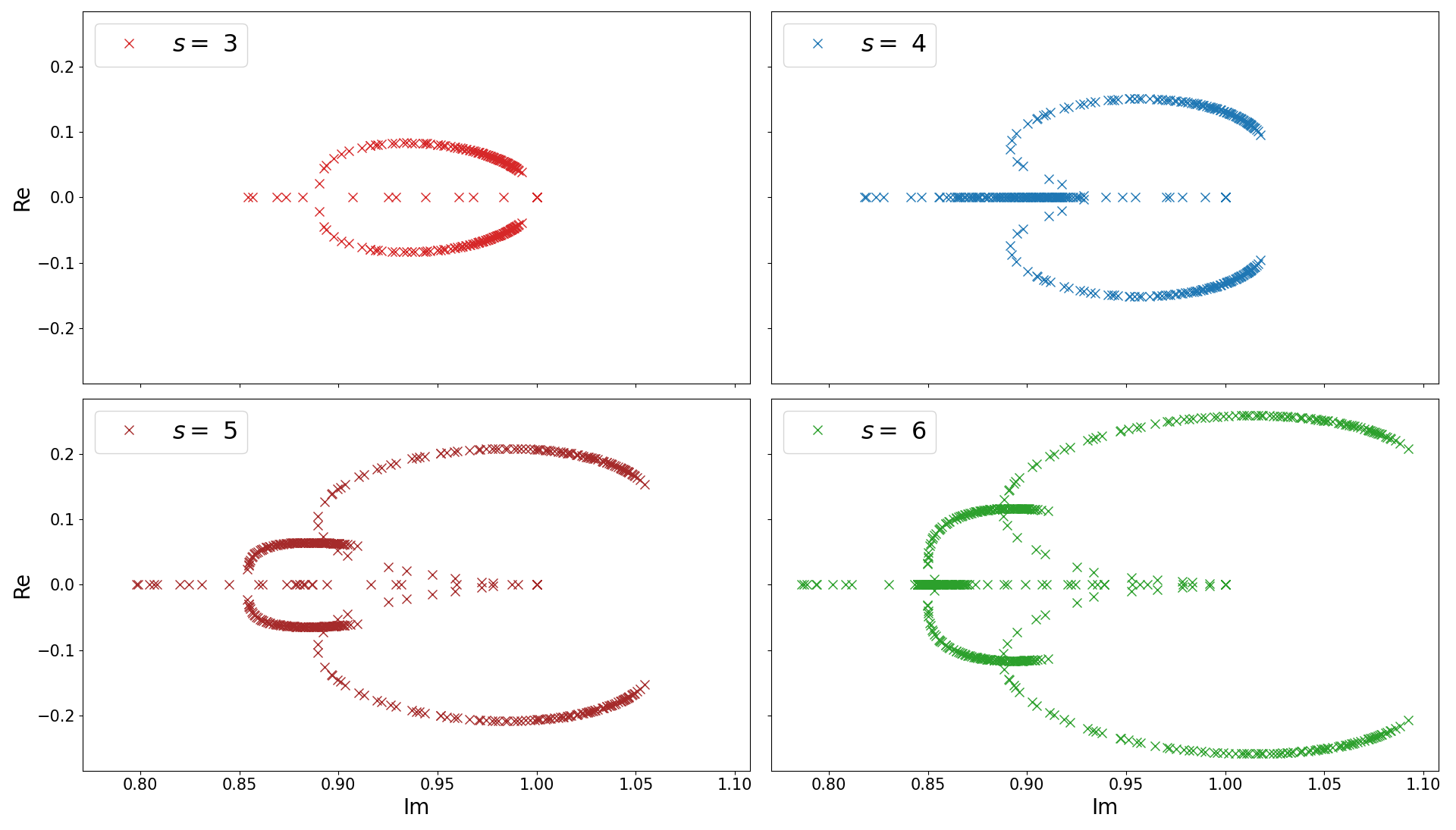}
\caption{The spectra of the preconditioned system $\mathcal{P}^{-1}\mathcal{A}$ for the preconditioner $\mathcal{P}$ given in~\eqref{eqn_secPolAprch_W2W1stencil_definition} below, using the RadauIIA IRK method. The spatial operator is the Laplacian on an irregular domain $\Omega$ with various boundary conditions (Dirichlet, Neumann and Robin), discretized using conforming P1 FEM, see~\cite[Section 4, Example 2]{outrata2023irksstage} for detailed description.}\label{fig_secIntro_GMRESspec}
\end{figure}

\noindent For the so-called fully IRK methods, the Butcher matrix $A = [a_{ij}]$ is dense, i.e., $a_{ij}\neq 0$ for all $i,j=1,\dotsc, s$, so that~\eqref{eqn_secIntro_StageEqntsInBasicForm} becomes a rather large system. We start by rewriting~\eqref{eqn_secIntro_StageEqntsInBasicForm} using a Kronecker product notation as

\begin{equation}\label{eqn_secIntro_IkronM_m_AkronK}
\left( I\otimes M + \tau A\otimes K \right) \mathbf{k}^{(m)} = \left( I\otimes M\right) \mathbf{k}^{(m-1)} + \begin{bmatrix} \mathbf{b}^{(ST)}(t_{m-1} + c_1\tau)^T,\dotsc , \mathbf{b}^{(ST)}(t_{m-1} + c_s\tau)^T \end{bmatrix}^T,
\end{equation}

\noindent where $\mathbf{k}^{(m-1)} = \begin{bmatrix} (\mathbf{k}^{(m-1)}_1)^T,\cdots , (\mathbf{k}^{(m-1)}_s)^T \end{bmatrix}^T$. Next, we transform~\eqref{eqn_secIntro_IkronM_m_AkronK} by factoring out $A\otimes I$ on the left-hand side to the left and multiplying the equation with its inverse, obtaining
\begin{equation}\label{eqn_secIntro_AinvkronM_m_IkronK}
\underbrace{ \left( A^{-1}\otimes M + \tau I \otimes K \right) }_{=: \, \mathcal{A}} \mathbf{k}^{(m)} = \mathbf{b}^{(m-1)}_{\mathrm{IRK}},
\end{equation}
\noindent so that $\mathbf{b}^{(m-1)}_{\mathrm{IRK}}$ is defined as the right-hand side vector in~\eqref{eqn_secIntro_IkronM_m_AkronK} multiplied by $A^{-1}\otimes I$ from the left. To the best of our knowledge this was first proposed by Butcher in~\cite{butcher1976implementation} and later became the standard, see~\cite{Howle2021new,howle2022efficient,Neytcheva2020,southworth2022fastI,southworth2022fastII,
dravins2024stageparallel_1,dravins2024stageparallel_2,dravins2024spectral,pazner2017stage}; the intuition on why this is useful is summarized in, e.g.,~\cite[Section 4]{dravins2024stageparallel_1} or~\cite[Section 4.1]{pazner2017stage}. The go-to solver for~\eqref{eqn_secIntro_AinvkronM_m_IkronK} has been the GMRES method (see~\cite[Sections 2.2, 5.7 and 5.10]{liesen2013krylov}) with a suitable preconditioner and in~\cite{staff2006preconditioning}, the authors considered a generic class of block-preconditioners (there formulated for~\eqref{eqn_secIntro_IkronM_m_AkronK}) given by
\begin{equation}\label{eqn_secIntro_PreconditionerWithAtilde}
\mathcal{P} = \tilde{A}\otimes M + \tau I \otimes K. 
\end{equation}
\noindent The authors proposed some convenient choice of $\tilde{A}$, which helped to stimulate the development of related preconditioners, notably the works~\cite{Howle2021new,howle2022efficient,Neytcheva2020,dravins2024stageparallel_2}. As the GMRES convergence is in practice often linked with the spectral properties of the (preconditioned) system matrix (see~\cite[Sections 2.2, 5.7 and 5.10]{liesen2013krylov}, also for the limits of this analysis), the authors presented also plots showing ``favorable properties''\footnote{In this context we would like to recall the classical result from~\cite{Greenbaum1996} that states that \emph{any} GMRES convergence can be observed for a system with a matrix with a given spectrum, i.e., spectrum on its own is \emph{not} sufficient to say anything about GMRES behavior. In practice, however, this is \emph{rarely} observed and many GMRES users use the folklore of ``better-clustered spectrum suggests faster convergence'', without further, case-specific justification. More details for this particular setting can be found in~\cite{outrata2023irksstage}.} of the spectrum. For illustration, we include an example in Figure~\ref{fig_secIntro_GMRESspec}.

%%%%%%%%%%%%%%%%%%%%%%%%%%%%%%%%%%%%%%%%%%%%%%%%%%%%%%%%%%%%
%%%%%%%%%%%%%%%%%%%%%%%%%%%%%%%%%%%%%%%%%%%%%%%%%%%%%%%%%%%%

As far as we are aware, there have been two independent series of works that aimed at analyzing preconditioners of the type~\eqref{eqn_secIntro_PreconditionerWithAtilde} -- namely their spectral properties -- the first coming from the group originally led by the late Owe Axelsson (\cite{dravins2024stageparallel_1,dravins2024spectral}) and the second from the group of Martin J. Gander (\cite{outrata2022IRKILAS,outrata2023irksstage}). The purpose of this work is to analyze their overlap and show their equivalence.

Both of these works are based on the spectral properties of the matrix pencil $\tau K-\mu M$ (we replace the standard symbol $\lambda$ for the generalized eigenvalue by $\mu$), i.e., on the eigendecomposition of the matrix $\tau M^{-1}K$. As the stiffness and mass matrices come from a discretization of a coercive, self-adjoint, bounded operator $\mathcal{L}$, we assume that the pencil $\tau K-\mu M$ is symmetric and positive-definite so that there exist an $M$-orthogonal eigenbasis $\mathbf{q}_1,\dotsc ,\mathbf{q}_n$ of $\tau K$, i.e.,
\begin{equation}\label{eqn_secIntro_EigenDecompOf_tauMinvK}
M^{-1}K = Q \begin{bmatrix}\mu_1 \\ & \ddots \\ && \mu_n \end{bmatrix} Q^T =:  QDQ^T ,
\quad \mathrm{with} \quad Q = [\mathbf{q}_1,\dotsc ,\mathbf{q}_n], \; D = \mathrm{diag}(\mu_1,\dotsc , \mu_n),
\end{equation}
\noindent where $0<\mu_1\leq \dotsc \leq \mu_n$ are the generalized eigenvalues of the pencil $\tau K-\mu M$, i.e., of the matrix $\tau M^{-1}K$, see~\cite[Sections 2.3 and 5]{bai2000templates} for more details and further references.

\section{Formulas for the eigenproperties of the preconditioned system}

\subsection{Polynomial approach}\label{sec_PolAprch}
This direction was initiated by the following observation of prof. Axelsson\,~\newline

\begin{center}
\begin{minipage}{.8\textwidth}
\textit{Taking the RadauIIA Butcher matrix $A$ for any $s=2,\dotsc ,10$, the lower-triangular part of $A^{-1}$ is dominating the rest of the matrix (in terms of magnitude of the entries). Similarly, the $LD$ factor of the LDU factorization of $A^{-1}$ is dominating the $U$ factor (in both the spectral and the Frobenius norm).}
\end{minipage}
\end{center}
\,~\newline

\noindent Heuristically, this suggests that these carry the majority of the information of $A^{-1}$ for the RadauIIA method and as such are reasonable choices for $\tilde{A}$ for the preconditioner construction. This has been numerically observed also in~\cite{Howle2021new,howle2022efficient}, only for $A$ itself, rather than $A^{-1}$. Writing the LDU factorization $A^{-1} = \tilde{L}\tilde{D}U$ we set $L:=\tilde{L}\tilde{D} = [l_{ij}]$ so that $U$ is a unit upper-triangular matrix, i.e., we have $ U = I + \hat{U}$ for some strictly upper-triangular matrix $\hat{U} = [\hat{u}_{ij}]$. Taking $\tilde{A} = L$ in~\eqref{eqn_secIntro_PreconditionerWithAtilde}, we obtain
\begin{equation}\label{eqn_secPolAprch_W2W1stencil_definition}
\begin{aligned}
\mathcal{P}^{-1}\mathcal{A}
&= \left( L\otimes M + \tau I \otimes K \right)^{-1} \left( L\otimes M + I \otimes \tau K  + (L\hat{U})\otimes M \right) \\
&= I + \left( I  + \tau L^{-1} \otimes M^{-1}K \right)^{-1} \left( \hat{U}\otimes I \right) =: I + W_1^{-1}W_2, \\
\end{aligned}
\end{equation}
\noindent see~\cite[equation (14)]{dravins2024stageparallel_1} or~\cite[Section 4]{dravins2024spectral} for further details. In \cite{dravins2024stageparallel_1}, the authors consider the RadauIIA IRK and show that for $s=2$ the eigenvalues of $\mathcal{P}^{-1}\mathcal{A}$ lie inside the disc centered at $1$ with the diameter $\|\hat{U}\|<1$ and conjecture this localization for general $s$. In~\cite[Section 4]{dravins2024spectral}, the authors improve the localization of the eigenvalues of $\mathcal{P}^{-1}\mathcal{A}$ for the RadauIIA IRK, using the matrix pencil $W_2-\lambda W_1$ (notice that $W_1,W_2$ are neither symmetric nor positive-definite). Namely, having a generalized eigenpair $(\lambda,\mathbf{v})$ of the pencil $W_2-\lambda W_1$ we have
\begin{equation}\label{eqn_secPolAprch_PinvA_W2W1_SpctrChar}
\mathcal{P}^{-1}\mathcal{A} \mathbf{v} = (1+\lambda) \mathbf{v} 
\quad \iff \quad
W_2 \mathbf{v} = \lambda W_1 \mathbf{v}.
\end{equation}
\noindent The authors then work with $W_2 \mathbf{v} = \lambda W_1 \mathbf{v}$ as with an equation for $(\lambda,\mathbf{v})$, i.e., the aim is to solve 
\begin{equation}\label{eqn_secPolAprch_W2W1stencil_GEformulation}
\begin{gathered}
\begin{bmatrix}
0 & \hat{u}_{1,2}I & \hdots & \hat{u}_{1,s}I \\
& \ddots & \ddots & \vdots \\
&& \ddots & \hat{u}_{s-1,s}I \\
&&& 0
\end{bmatrix} \begin{bmatrix} \mathbf{v}_1\\ \vdots \\ \vdots \\ \mathbf{v}_s \end{bmatrix} = \lambda \left( \begin{bmatrix} \mathbf{v}_1\\ \vdots \\ \vdots \\ \mathbf{v}_s \end{bmatrix} + \begin{bmatrix}
\ell_{11}\cdot \tau M^{-1}K \\
\vdots & \ddots \\
\vdots && \ddots \\
\ell_{s1}\cdot \tau M^{-1}K & \hdots  & \hdots & \ell_{ss}\cdot \tau M^{-1}K
\end{bmatrix} \begin{bmatrix} \mathbf{v}_1\\ \vdots \\ \vdots \\ \mathbf{v}_s \end{bmatrix} \right)
\end{gathered}
\end{equation}
\noindent for $(\lambda,\mathbf{v})$. In~\cite[Section 4.3]{dravins2024spectral}, the authors first observe $n$ ``trivial'' solutions corresponding to $\lambda=0$. To resolve the remaining $(s-1)n$ eigenpairs, the authors proceed with a \emph{symbolic} block backward substitution. This is presented for $s=2,3$ in Sections~4.1 and~4.2 and by analogy the authors argue that such a routine can be carried out for any $s$. As there is little space devoted to this argument in~\cite[Section 4.3]{dravins2024spectral}, we derive it here ourselves, adhering to the techniques used in~\cite[Section 4]{dravins2024spectral}. 

Looking at the $i$-th ($1\leq i \leq s$) block-row in~\eqref{eqn_secPolAprch_W2W1stencil_GEformulation}, we can rearrange it as
\begin{equation}\label{eqn_secPolAprch_vi_EqnFromGEreshuffled}
\lambda(I+\ell_{ii}\cdot \tau M^{-1}K) \mathbf{v}_i  = \sum\limits_{j=i+1}^s \hat{u}_{i,j} \mathbf{v}_j - \sum\limits_{j=1}^{i-1} \lambda \ell_{i,j} \cdot \tau M^{-1}K \mathbf{v}_j,
\end{equation}
\noindent and notice that for $i=s$, $\lambda$ factors out from both sides so that $\mathbf{v}_s$ can be expressed only in terms of $\mathbf{v}_1,\dotsc ,\mathbf{v}_{s-1}$, i.e., independent of $\lambda$. This corresponds to the first  step of the mentioned symbolic backward substitution. We proceed with the substitution for $i=s-1,\dotsc ,1$. When dealing with the $i$-th block-row we assume to have expressions for $\mathbf{v}_{s},\dotsc ,\mathbf{v}_{i+1}$ in terms of $\mathbf{v}_1,\dotsc ,\mathbf{v}_{i}$ and we insert these into~\eqref{eqn_secPolAprch_vi_EqnFromGEreshuffled}, obtaining 
\begin{equation}\label{eqn_secPolAprch_vi_EqnFromGEreshuffled_withGiHi}
\left[ \lambda(I+\ell_{ii}\cdot \tau M^{-1}K) - H^{(i)} \right] \mathbf{v}_i =
\sum\limits_{j=1}^{i-1} \left( G^{(i)}_j - \lambda \ell_{i,j} \cdot \tau M^{-1}K \right) \mathbf{v}_j,
\end{equation}
\noindent for some appropriate $n$-by-$n$ matrices $H^{(i)}, G^{(i)}_j$ based on~\eqref{eqn_secPolAprch_vi_EqnFromGEreshuffled}. For example, 
\begin{equation}\label{eqn_secPolAprch_v1_examplHsm1Gsm1}
\begin{aligned}
H^{(s-1)} &= -\hat{u}_{s-1,s} (I+\ell_{ss}\cdot \tau M^{-1}K)^{-1} \ell_{s,s-1}\cdot \tau M^{-1}K,\\
G^{(s-1)}_j &= -\hat{u}_{s-1,s} (I+\ell_{ss}\cdot \tau M^{-1}K)^{-1} \ell_{s,j}\cdot \tau M^{-1}K.
\end{aligned}
\end{equation}
\noindent From~\eqref{eqn_secPolAprch_vi_EqnFromGEreshuffled_withGiHi} we obtain an expression for $\mathbf{v}_i$ in terms of $\mathbf{v}_{i-1},\dotsc ,\mathbf{v}_1$ and after $s$ steps arrive at
\begin{equation}\label{eqn_secPolAprch_v1_EqnFtildev1eq0}
\underbrace{ \, \left[ \lambda(I+\ell_{11}\cdot \tau M^{-1}K) - H^{(1)} \right] \,}_{=: \; \tilde{F}_{\tau M^{-1}K}(\lambda)} \mathbf{v}_1 = 0.
\end{equation}
\noindent This process ``sandwiches'' the inversion of the matrices $\lambda(I+\ell_{ii}\cdot \tau M^{-1}K) - H^{(i)}$ and the multiplication with matrices $G^{(i)}_j - \lambda \ell_{i,j} \cdot \tau M^{-1}K$. In other words, $\tilde{F}_{\tau M^{-1}K}(\lambda)$ is built from the building blocks of matrices $\tau M^{-1}K$ and $I$ by repeated applications of linear combination and inversion. This is, at least in our understanding, the key observation behind the analysis in~\cite{dravins2024spectral} and has several relevant implications.

First, in this sense $\tilde{F}_{\tau M^{-1}K}(\lambda)$ is a rational function of $\tau M^{-1}K$ and also a rational function of $\lambda$. Second, all matrix products in the definition of $\tilde{F}_{\tau M^{-1}K}(\lambda)$ commute, for example we have
\begin{equation}\label{eqn_secPolAprch_v1_examplHsm2}
\begin{aligned}
H^{(s-2)} &= \hat{u}_{s-2,s} G^{(s-1)}_{s-2} - \left( \hat{u}_{s-2,s-1} + \hat{u}_{s-2,s}  (I+\ell_{ss}\cdot \tau M^{-1}K)^{-1} \ell_{s,s-1}\cdot \tau M^{-1}K \right) \cdot \\ 
			& \; \cdot \left[ \lambda(I+\ell_{s-1,s-1}\cdot \tau M^{-1}K) - H^{(s-1)} \right]^{-1} \left( G^{(s-1)}_{s-2} - \lambda \ell_{s-1,s-2} \cdot \tau M^{-1}K \right) \\
%%%%%%%%%%%%%%%%%
&= \hat{u}_{s-2,s} G^{(s-1)}_{s-2} -  \left[ \lambda(I+\ell_{s-1,s-1}\cdot \tau M^{-1}K) - H^{(s-1)} \right]^{-1} \cdot \\ 
			& \;  \cdot \left( \hat{u}_{s-2,s-1} + \hat{u}_{s-2,s}  (I+\ell_{ss}\cdot \tau M^{-1}K)^{-1} \ell_{s,s-1}\cdot \tau M^{-1}K \right) \left( G^{(s-1)}_{s-2} - \lambda \ell_{s-1,s-2} \cdot \tau M^{-1}K \right).
\end{aligned}
\end{equation}
\noindent In particular, the inversed matrices ``sandwiched inside'' $H^{(1)}$ can freely ``travel to the left'' within the sandwiches, i.e., we can move all the inversed matrices in the matrix products inside $\tilde{F}_{\tau M^{-1}K}(\lambda)$ (i.e., inside $H^{(1)}$) ``to the left'', similarly to~\eqref{eqn_secPolAprch_v1_examplHsm2}. After this rearrangement is finished, we can multiply (also from the left) with those matrices, transforming~\eqref{eqn_secPolAprch_v1_EqnFtildev1eq0} so as to get rid of all inversions. For example, in~\eqref{eqn_secPolAprch_v1_examplHsm2} we would first multiply by $\lambda(I+\ell_{s-1,s-1}\cdot \tau M^{-1}K) - H^{(s-1)}$ and then with $I+\ell_{ss}\cdot \tau M^{-1}K$ (which is also present inside $H^{(s-1)}$). This corresponds to expanding back the backward substitution process and in this process we transform~\eqref{eqn_secPolAprch_v1_EqnFtildev1eq0} to a new equation, let us denote it
\begin{equation}\label{eqn_secPolAprch_v1_EqnFv1eq0}
F_{\tau M^{-1}K}(\lambda) \mathbf{v}_1 = 0.
\end{equation}
\noindent Already for $s=3$ (with $H^{(1)} = H^{(s-2)}$ given in~\eqref{eqn_secPolAprch_v1_examplHsm2}) the calculation becomes somewhat tedious and it is a downright unpleasant chore for larger number of stages. However, it allows us to focus on the quantity of interest -- the matrix $F_{\tau M^{-1}K}(\lambda)$. 

By construction, each of the inversed matrices consisted of linear combination of (a) terms linear in $\lambda$ and (b) other inversed matrices. As we did $s$-step backward substitution, the above process of transforming $\tilde{F}_{\tau M^{-1}K}(\lambda)$ into $F_{\tau M^{-1}K}(\lambda)$ introduces at most $s$ multiplications -- the first $s-1$ of them will introduce a linear factor in $\lambda$, while the last (corresponding to the bottom-most block-row of~\eqref{eqn_secPolAprch_W2W1stencil_GEformulation}) is independent of $\lambda$. In other words, $F_{\tau M^{-1}K}(\lambda)$ is a polynomial function in $\lambda$ of degree $s-1$. In fact, the construction showcases that $F_{\tau M^{-1}K}(\lambda)$ is the numerator of the rational function $\tilde{F}_{\tau M^{-1}K}(\lambda)$ with respect to $\lambda$. Also by construction,  
\begin{equation}\label{eqn_secPolAprch_KerFtild_eq_KerF}
\mathrm{Ker}\left( F_{\tau M^{-1}K}(\lambda) \right) = \mathrm{Ker}\left( \tilde{F}_{\tau M^{-1}K}(\lambda) \right),
\end{equation}
\noindent since we always multiplied by non-singular matrices with trivial kernels. Linking this back to the preconditioned system, we observe that $(1+\lambda,\mathbf{v})$ is an eigenpair of the preconditioned system $\mathcal{P}^{-1}\mathcal{A}$ if and only if $\mathbf{v}_1 \in \mathrm{Ker}\left( F_{\tau M^{-1}K}(\lambda) \right)$. This shows there has to be a very close link between the matrix polynomial $F_{\tau M^{-1}K}(\lambda)$ and the characteristic polynomial of $\mathcal{P}^{-1}\mathcal{A}$ (after the simple change of variables $\lambda = 1+\lambda$, see~\eqref{eqn_secPolAprch_PinvA_W2W1_SpctrChar}).

A third consequence of this special structure of $\tilde{F}_{\tau M^{-1}K}(\lambda)$ (which extends to $F_{\tau M^{-1}K}(\lambda)$) is that $\tilde{F}_{\tau M^{-1}K}(\lambda)$ must diagonalize in the eigenbasis $Q$ of $\tau M^{-1}K$, i.e., we can rewrite~\eqref{eqn_secPolAprch_v1_EqnFv1eq0} as
\begin{equation*}
\begin{aligned}
Q^T \tilde{F}_{\tau M^{-1}K}(\lambda) Q Q^T\mathbf{v}_1 \equiv \begin{bmatrix} \tilde{F}_{\mu_1}(\lambda)\\ & \ddots \\ && \tilde{F}_{\mu_n}(\lambda) \end{bmatrix} Q^T\mathbf{v}_1 &= 0, \\
Q^T F_{\tau M^{-1}K}(\lambda) Q Q^T\mathbf{v}_1 \equiv \begin{bmatrix} F_{\mu_1}(\lambda)\\ & \ddots \\ && F_{\mu_n}(\lambda) \end{bmatrix} Q^T\mathbf{v}_1 &= 0,
\end{aligned}
\end{equation*}
\noindent where $\tilde{F}_{\mu_k}(\lambda)$ (or $F_{\mu_k}(\lambda)$) are now \emph{scalar} rational (or polynomial) functions of $\lambda$ with the precise structure of $\tilde{F}_{\tau M^{-1}K}(\lambda)$ (or $F_{\tau M^{-1}K}(\lambda)$), only replacing the matrices $\tau M^{-1}K$ and $I$ by $\mu_k$ and $1$. Hence, the coefficients of $\tilde{F}_{\mu_k}(\lambda)$ (or $F_{\mu_k}(\lambda)$) are themselves rational functions in $\mu_k$. Crucially, this showcases the connection between the matrix $F_{\tau M^{-1}K}(\lambda)$ and the characteristic polynomial of $\mathcal{P}^{-1}\mathcal{A}$. By the same argument as in~\eqref{eqn_secPolAprch_KerFtild_eq_KerF} we get
\begin{equation*}
\mathrm{det}\left( F_{\tau M^{-1}K}(\lambda) \right) = \alpha_{\tau M^{-1}K} \mathrm{det}\left( \tilde{F}_{\tau M^{-1}K}(\lambda) \right),
\end{equation*}
\noindent for some scaling coefficient $\alpha_{\tau M^{-1}K}\neq 0$ and thereby we see that if $\mathrm{det}\left( F_{\tau M^{-1}K}(\lambda) \right) = 0$, then $\lambda$ is one of the $s(n-1)$ generalized eigenvalues of the pencil $W_2-\lambda W_1$ we are looking for. Writing
\begin{equation}\label{eqn_secPolAprch_CharPol_W1W2stencil}
\hat{p}_{\mathrm{char}}(\lambda) := \det \left( F_{\tau M^{-1}K}(\lambda) \right) =
\det \begin{bmatrix} F_{\mu_1}(\lambda)\\ & \ddots \\ && F_{\mu_n}(\lambda) \end{bmatrix} =
\prod\limits_{k=1}^{n} F_{\mu_k}(\lambda),
\end{equation}
\noindent it follows that $p_{\mathrm{char}} (\lambda) : = \lambda^n \hat{p}_{\mathrm{char}} (\lambda)$ is (up to a rescaling) the characteristic polynomial of the pencil $W_2-\lambda W_1$, and therefore $p_{\mathrm{char}} (\lambda-1)$ is (up to a rescaling) the characteristic polynomial of the preconditioned system $\mathcal{P}^{-1}\mathcal{A}$. Due to the structure of $\hat{p}_{\mathrm{char}}(\lambda)$ in~\eqref{eqn_secPolAprch_CharPol_W1W2stencil} it becomes natural to index the eigenvalues of $\mathcal{P}^{-1}\mathcal{A}$ (i.e., the zeros of $p_{\mathrm{char}}$) by the eigenvalues $\mu_k$ of $\tau M^{-1}K$.

\begin{remark}
While not stated this way, the above was clearly grasped in~\cite[Sections 4 and 5]{dravins2024spectral}. There, the above derivation is essentailly skipped by taking $\mathbf{v}_1 = \mathbf{q}_1$. This somewhat simplifies the calculation compared to~(\ref{eqn_secPolAprch_vi_EqnFromGEreshuffled_withGiHi}-\ref{eqn_secPolAprch_v1_EqnFv1eq0}) but not to an agreeable level. This is also illustrated by the fact that the numerical experiments in~\cite{dravins2024spectral} are carried out only for $s=2,3$. Comparing with~\cite[Sections 4.3 and Theorem 5.1]{dravins2024spectral}, the above is, in our eyes, slightly more constructive way to arrive at the same result. We also believe that it highlights more clearly the key mechanic of the derivation.
\end{remark}

However, since the analysis in~\cite[Sections 4 and 5]{dravins2024spectral} relies on the  \emph{construction} of $p_{\mathrm{char}}$, the above seems not fully satisfactory -- as mentioned above, it is a laborious task to obtain the symbolic formulas beyond $s=2,3$ and this practical aspect was not addressed; authors simply state that this symbolic approach must work for any $s$.

A slight simplification can be achieved as follows. We return to~\eqref{eqn_secPolAprch_W2W1stencil_GEformulation} and pass \emph{blockwise} into the eigenbasis $Q$ directly there, obtaining
\begin{equation}\label{eqn_secPolAprch_W2W1Stencil_GEformulation_KronBfrSclrVersion}
\left( \hat{U}\otimes I_n \right) \mathbf{w} = \lambda \left( I_s\otimes I_n + L^{-1} \otimes D \right) \mathbf{w},
\end{equation}
\noindent with $\mathbf{w}=(I\otimes Q^T)\mathbf{v}$. This transposes the block backward substitution of~\eqref{eqn_secPolAprch_W2W1stencil_GEformulation} into $n$ independent scalar ones,
\begin{equation}\label{eqn_secPolAprch_W2W1Stencil_GEformulation_SclrVersion}
\hat{U} \mathbf{s}_{\mu_k} = \lambda \left( I_s + \mu_k L^{-1} \right) \mathbf{s}_{\mu_k},
\end{equation}
\noindent again, parametrized by $\mu_k, k = 1,\dotsc,n$. By construction, the characteristic polynomial of the matrix pencil $\hat{U} - \lambda \left( I_s + \mu_k L^{-1} \right)$ is (up to a rescaling) identical to $F_{\mu_k}(\lambda)$, i.e., we can calculate the coefficients of $F_{\mu_k}(\lambda)$ using a scalar version of the symbolic backward substitution -- arriving at formulas from~\cite[Sections 4.2]{dravins2024spectral} for the case $s=3$ \emph{without} the blockwise backward substitution in~\eqref{eqn_secPolAprch_W2W1stencil_GEformulation}.

Obtaining the spectrum of $\mathcal{P}^{-1}\mathcal{A}$ (or an estimate of it) then reduces to calculating (or estimating) $\mu_1,\dotsc, \mu_n$, then evaluating the symbolic formulas to get the coefficients $\mathbf{c}_1,\dotsc,\mathbf{c}_n \in \mathbb{R}^{s}$ of $F_{\mu_1}(\lambda),\dotsc ,F_{\mu_n}(\lambda)$ and calculating the roots of these polynomials using a standard numerical software, e.g., $\mathtt{numpy.roots}(\mathbf{c}_k)$. 

A direct observation is that, from the numerical perspective, finding the roots of a higher-degree polynomial written in the monomial basis can be a poorly-conditioned problem with respect to the polynomial coefficients. This is a well-known problem and algorithms for approximation of the roots of a polynomial given by a coefficient vector $\mathbf{c} \in \mathbb{R}^d$, such as $\mathtt{numpy.roots}()$, circumvent the issue by calculating eigenvalues of the so-called companion matrix $\mathbf{C}_{\mathbf{c}} \in \mathbb{R}^{d\times d}$. In particular, the costs of finding the roots of a polynomial with coefficients $\mathbf{c} \in \mathbb{R}^{d}$ are governed by the costs of solving a $d$-dimensional eigenvalue problem. In our setting, this means that the costs of finding the eigenvalues $\lambda^{(k)}_1,\dotsc, \lambda^{(k)}_{s-1}, k=1,\dotsc ,n$ of the matrix pencil $W_2-\lambda W_1$ (or their estimates) after all the symbolic work has been done still amounts to solving $n$ independent eigenvalue problems of size $(s-1)$-by-$(s-1)$. Taking a step back, we notice that we have already expressed $\lambda^{(k)}_1,\dotsc, \lambda^{(k)}_{s-1}$ as the eigenvalues of $n$ independent $s$-by-$s$ generalized eigenvalue problems -- in the equation~\eqref{eqn_secPolAprch_W2W1Stencil_GEformulation_SclrVersion}. Moreover, their construction is \emph{direct}, given the values (or approximations of) $\mu_1,\dotsc ,\mu_n$ -- no symbolic calculations are needed.

In other words, what we gain by the explicit construction of $F_{\mu_1}(\lambda),\dotsc ,F_{\mu_n}(\lambda)$ compared to numerically solving~\eqref{eqn_secPolAprch_W2W1Stencil_GEformulation_SclrVersion} is the reduction of the size of each of the $n$ independent eigenvalue problems by one and the fact that we can use eigenvalue solvers tailored to companion matrices rather than generic generalized eigenvalue solver. Given that these problems are \emph{very} small, we believe it to be more efficient and numerically stable, to avoid the explicit construction of $F_{\mu_1}(\lambda),\dotsc ,F_{\mu_n}(\lambda)$ and instead treat these only implicitly by using an eigenvalue solver for~\eqref{eqn_secPolAprch_W2W1Stencil_GEformulation_SclrVersion}. Crucially, this is not present in~\cite{dravins2024spectral} and it is the keystone between the analysis above and the one presented in~\cite{outrata2023irksstage} as we shall see next.

\subsection{Matrix approach}
This direction was initiated by prof. Gander based on the plenary talk of prof. Howle at the PRECOND conference in 2019 in Minneapolis (later, some of the results appeared in~\cite{Howle2021new}). Her group has been interested in preconditioners for the systems~\eqref{eqn_secIntro_StageEqntsInBasicForm} for problems like~\eqref{eqn_secIntro_ParabPDE} and she presented some strong numerical results for the block preconditioners of the type~\eqref{eqn_secIntro_PreconditionerWithAtilde} (and since has expanded the focus also to hyperbolic problems, see~\cite{howle2022efficient,outrata2025irkAman}). As our analysis was framed mainly by the setup in~\cite{Howle2021new}, we dealt with the preconditioned system
\begin{equation*}
\left( I\otimes M + \tau \tilde{A} \otimes K \right)^{-1}\left( I\otimes M + \tau A \otimes K \right) \mathbf{k}^{(m+1)} = \tilde{\mathbf{b}}^{(m)},
\end{equation*}
\noindent rather than~\eqref{eqn_secIntro_AinvkronM_m_IkronK} with the preconditioner~\eqref{eqn_secIntro_PreconditionerWithAtilde} -- here we will write the results in the ``notation'' of~\eqref{eqn_secIntro_AinvkronM_m_IkronK} and~\eqref{eqn_secIntro_PreconditionerWithAtilde}. 

Similarly to~\cite{dravins2024stageparallel_1,dravins2024spectral}, we also first studied the case $s=2$ in detail in~\cite{outrata2022IRKILAS} and only later generalized the results for arbitrary $s$ in~\cite{outrata2023irksstage}, see also~\cite{outrata2022phdthesis}.

The analysis is built using the same observation that allowed us to pass from~\eqref{eqn_secPolAprch_W2W1Stencil_GEformulation_KronBfrSclrVersion} to~\eqref{eqn_secPolAprch_W2W1Stencil_GEformulation_SclrVersion}, which is based on the Kronecker product properties. In particular, we can write
\begin{equation*}
\begin{aligned}
\mathcal{P}^{-1}\mathcal{A}
&= \left( L\otimes M + \tau I \otimes K \right)^{-1} \left( A^{-1}\otimes M + I \otimes \tau K\right) \\
&= (I\otimes Q^T)^{-1} \underbrace{\, \left( L\otimes I  + I\otimes D \right)^{-1} \left( A^{-1}\otimes I + I\otimes D \right) \,}_{=: \; X_{\tau M^{-1}K}} (I\otimes Q^T). \\
\end{aligned}
\end{equation*}
\noindent By construction, the matrix $X_{\tau M^{-1}K}$ is an $ns$-by-$ns$ block matrix, each block $X_{\tau M^{-1}K}^{(ij)} \in \mathbb{R}^{n\times n}$ being \emph{a diagonal matrix}\footnote{We choose the subscript $\tau M^{-1}K$ by analogy to Section~\ref{sec_PolAprch}, as they play a similar role in the respective expositions.}, $X_{\tau M^{-1}K}^{(ij)} = \mathrm{diag}(x^{(ij)}_1,\dotsc ,x^{(ij)}_n)$. This is a very special sparsity structure and can be transformed into one we are perhaps more used to -- we can permute  $X_{\tau M^{-1}K}$ into a block-diagonal matrix. That is, there exists a permutation matrix $\Pi$ such that
\begin{equation*}
\Pi^T X_{\tau M^{-1}K} \Pi \equiv \Pi^T 
\begin{bmatrix}
X_{\tau M^{-1}K}^{(11)} & \hdots & X_{\tau M^{-1}K}^{(1s)} \\
\vdots & \ddots & \vdots\\
X_{\tau M^{-1}K}^{(s1)} & \hdots & X_{\tau M^{-1}K}^{(ss)}
\end{bmatrix} \Pi = 
\begin{bmatrix}
X_{1} \\ & \ddots \\ && X_{n}
\end{bmatrix},
\end{equation*}
\noindent where $X_{1}, \dotsc ,X_{n}$ are given as
\begin{equation}\label{eqn_secMtrxAprch_Xk_def_eq_sBYsForm}
X_{k} = \begin{bmatrix} x^{(11)}_k & \hdots & x^{(s1)}_k \\ \vdots & \ddots & \vdots \\ x^{(s1)}_k & \hdots & x^{(ss)}_k \end{bmatrix} = (L+\mu_k)^{-1}(A^{-1}+\mu_k) \in \mathbb{R}^{s\times s},
\end{equation}
\noindent for $k=1,\dotsc , n$, see~\cite[Lemma 3.1]{outrata2023irksstage}. Importantly, $X_k$ contains precisely the terms in $X_{\tau M^{-1}K}$ that depend on the $k$-th eigenspace corresponding to $\mu_k$ and is independent of $\mu_1,\dotsc ,\mu_{k-1},\mu_{k
+1}, \dotsc,\mu_{n}$. Hence, we will use the notation $X_{\mu_k}$ rather than  $X_{k}$, similarly to Section~\ref{sec_PolAprch}. By construction, having an eigenpair $(\lambda,\mathbf{s})$ of $X_{\mu_k}$ we have that
\begin{equation}\label{eqn_secMtrxAprch_PinvA_EigenpairsBasedOn_Xmuk}
\mathcal{P}^{-1}\mathcal{A} (I\otimes Q^T)^{-1} \Pi (\mathbf{e}_k\otimes \mathbf{s}) = 
(I\otimes Q^T)^{-1} \Pi (\mathbf{e}_k\otimes X_{\mu_k}\mathbf{s}) = 
\lambda (I\otimes Q^T)^{-1} \Pi (\mathbf{e}_k\otimes \mathbf{s}),
\end{equation}
\noindent i.e., the eigenpairs of the preconditioned system are fully characterized by those of $X_{\mu_k} \in \mathbb{R}^{s\times s}$.

\begin{remark}
We would like to highlight that the matrices $X_{\mu_k}$ combined with~\eqref{eqn_secMtrxAprch_PinvA_EigenpairsBasedOn_Xmuk} turn out to be practically useful also for estimation of the standard GMRES bound
\begin{equation}\label{eqn_secMtrxAprch_GMRES_idealGMRES_poly_in_eigenvals}
\frac{\| \mathbf{r}_{\ell}\|}{\| \mathbf{r}_0\|} \leq  \; \kappa (S) \; \min\limits_{
\substack{\varphi(0)=1 \\ \mathrm{deg}(\varphi)\leq \ell} }  \max\limits_{1\leq i \leq sn} |\varphi(\lambda_i)|,
\end{equation}
\noindent where $S$ is the matrix of eigenvectors of $\mathcal{P}^{-1}\mathcal{A}$ and $\kappa (S)$ is its condition number. As we can see in~\cite[Section~4]{outrata2023irksstage}, the resulting GMRES convergence estimates are very descriptive even for larger number of stages. Crucially, though, the original version of the work contained a simple error in the last conclusion drawn from~\eqref{eqn_secMtrxAprch_PinvA_EigenpairsBasedOn_Xmuk}, namely in evaluation of $\kappa (S)$.
The corrected version has been submitted since and is also present in the preprints (links are on the personal webpages of both of the authors).
\end{remark}

\subsection{Connecting the two}
The first point of contact between the two groups happened at the SIAM LA meeting in Paris in 2024, where I have met Ivo Dravins who was presenting the materials from~\cite{dravins2024spectral}. This lead to several useful discussions. Among other things, these discussions led me to write this text to clearly connect the two approaches for analyzing the spectrum of $\mathcal{P}^{-1}\mathcal{A}$ -- in fact, we will see next how to map one onto the other.

The key to seeing the two approaches as one is again the Kronecker product structure and arithmetic. Comparing~\eqref{eqn_secMtrxAprch_Xk_def_eq_sBYsForm} with~\eqref{eqn_secPolAprch_W2W1Stencil_GEformulation_SclrVersion}, we see two sets of $n$ eigenvalue problems of the size $s$-by-$s$, parametrized by $\mu_k$. A direct calculation gives
\begin{equation}\label{eqn_secMtrxAprch_Xk_to_W2W1Stencil}
\begin{aligned}
X_{\mu_k} &= (L+\mu_k)^{-1}(A^{-1}+\mu_k) = (L+\mu_k)^{-1}(L(I+\hat{U})+\mu_k)\\
&= I + (L+\mu_k)^{-1}(L\hat{U})) = I + (I+L^{-1}\mu_k)^{-1}\hat{U}),
\end{aligned}
\end{equation}
\noindent the analogue of the manipulations in~\eqref{eqn_secPolAprch_W2W1stencil_definition}, only carried out with the $s$-by-$s$ blocks on the diagonal, after passing into the block basis $(I\otimes Q^T)\Pi$. Following this analogy one step further, we see that the generalized eigenvalue problems in~\eqref{eqn_secPolAprch_W2W1Stencil_GEformulation_SclrVersion} are simply a redressing of the eigenvalue problems with $X_{\mu_k}$ in~\eqref{eqn_secMtrxAprch_Xk_def_eq_sBYsForm}. Relating this back to the objects featuring in~\cite{dravins2024spectral,outrata2023irksstage}, for any $k=1,\dotsc ,n$ we have
\begin{equation*}
\det \left( \lambda I - X_{\mu_k} \right) =  \det \left( (\lambda-1) I - (I+L^{-1}\mu_k)^{-1}\hat{U}\right).
\end{equation*}
\noindent Recalling that $F_{\mu_k}(\lambda)$ is obtained by factoring out $\lambda$ in the first step of the (block) backward substitution, we conclude that $(\lambda - 1) F_{\mu_k}(\lambda-1)$ is the characteristic polynomial of $X_{\mu_k}$ (up to a rescaling). In other words, the approach in~\cite{dravins2024spectral} led to assembling the characteristic polynomial of the matrices used in~\cite{outrata2023irksstage}, both posed in the block eigenbasis $I\otimes Q^T$.

In both of the works, the authors recognized that this approach is fully general with respect to the chosen IRK method, i.e., with respect to $A$, although the results naturally do depend on this choice. Moreover, the authors also recognized, that their respective formulations of the spectrum of $\mathcal{P}^{-1}\mathcal{A}$ reveal its additional structure -- the eigenvalues come from a single object -- the matrix $X_{\mu_k}$ (or the polynomial $p_{\mathrm{char}}$) -- parametrized by $\mu_k$. Since in many cases of interest the eigenvalues $\mu_k$ sample fairly densely some interval $(\mu_{\mathrm{min}},\mu_{\mathrm{max}})$, each eigenvalue (zero) of $X_{\mu_k}$ ($p_{\mathrm{char}}$) becomes a continuous function of $\mu_k$. In particular, the eigenvalues of $\mathcal{P}^{-1}\mathcal{A}$ necessarily form \emph{branches} -- in fact $s$ of them -- where each branch tracks out one of these ``smooth functions'' of $\mu_k$.

\begin{remark}
We can observe this in Figure~\ref{fig_secIntro_GMRESspec} -- the complex conjugate branches are densely sampled even for a very low mesh resolution ($h\approx 0.7$). For $s=3$ we see clearly two complex conjugate branches. There are also several real eigenvalues. In fact $N$ of them are exactly equal to $1+0i$ -- a theoretical result established in both~\cite{dravins2024spectral} and~\cite{outrata2023irksstage} -- and these constitute the third, branch (hence we indeed have $s$ branches). The real eigenvalues other than $1+0i$ highlight that, from a certain $\mu_k$ onwards, the matrices $X_{\mu_k}$ have real spectrum. In particular, there exists a $\hat{\mu}_k$ for which $X_{\hat{\mu}_k}$ has a single real eigenvalue $\hat{\lambda}_k$ with algebraic multiplicity equal to two (based on Figure~\ref{fig_secIntro_GMRESspec} we have $\hat{\lambda}_k \approx 0.89 + 0i$). In this way, Figure~\ref{fig_secIntro_GMRESspec} nicely illustrates that indeed the eigenvalues of $\mathcal{P}^{-1}\mathcal{A}$ are only continuous (and not, e.g., $\mathcal{C}^{1}$) functions of $\mu_k$. Similar behavior is present in all of the graphs in Figure~\ref{fig_secIntro_GMRESspec} -- we have $s$ branches, one of them degenerating into the single point $1+0i$ and the others appearing in complex conjugate pairs (if they are indeed complex); notably for $s=4,6$ there is an additional, purely real branch.
\end{remark}

\section{Conclusion and future work}
We have demonstrated how to transpose the analysis of each of the two groups into the language of the other, spelling out in detail how these seemingly different approaches map onto each other. However, in both~\cite{dravins2024spectral,outrata2023irksstage}, the spectral results are only a part of the story. In~\cite{dravins2024spectral} these are used together with the block locally Toeplitz (BLT) theory to establish interesting asymptotic results (which have been independently observed also in~\cite[Chapter 7]{outrata2022phdthesis}). In~\cite{outrata2023irksstage}, the focus is on using the obtained eigeninformation (or an estimate thereof) for estimating the GMRES behavior convergence of the preconditioned systems, using the potential theory and the Schwarz-Christoffel (SC) mapping in particular. We believe that these directions could and should be explored in conjunction as well -- the BLT theory can, in principle, provide estimates of the (fully complex) spectrum of the spatial operator, which for is the crucial piece needed for further generalization of the SC mapping approach for GMRES estimation for systems stemming from IRK applied to PDEs with more involved spatial operators. This is a work in progress and which we find quite enticing.

\section{Acknowledgment}
  The author would like to acknowledge the support of the PRIMUS grant \texttt{PRIMUS/25/SCI/022} and would like to thank Ivo Dravins for several inspiring discussions.

\bibliographystyle{plainnat}
\bibliography{MyBiblio}

\end{document}